\documentclass[letterpaper,10pt]{article}
\usepackage[utf8x]{inputenc}
\usepackage{amsthm,amsmath,amssymb,amscd}
\usepackage{graphicx}

\newtheorem{definicion}{Definition}
\newtheorem{teorema}{Theorem}

\renewcommand{\P}[1]{{\mathbf P}_#1}
\renewcommand{\S}[1]{{\mathbf S}_#1}

\newcommand{\R}{{\mathbb R}}
\newcommand{\N}{{\mathbb N}}

\title{Explicit Spectral Decimation for a Class of Self--Similar Fractals}
\author{Sergio A. Hern\'andez\and Federico Men\'endez--Conde}

\date{}

\begin{document}

\maketitle

\begin{center}

Centro de Investigaci\'on en Matem\'aticas 

Universidad Aut\'onoma del Estado de Hidalgo 

fmclara@uaeh.edu.mx

hasasaha@gmail.com
\end{center}

\begin{abstract}
The method of spectral decimation is applied to an infinite collection of self--similar fractals. The sets considered are a generalization of
the Sierpinski Gasket to higher dimensions; they belong to the class of nested fractals, 
and are thus very symmetric. An explicit construction is given to obtain formulas for the eigenvalues of the Laplace operator acting on these fractals.  
\end{abstract}
\medskip 

In 1989, J. Kigami \cite{Kigami89} gave an analytic definition of a Laplace operator acting on the Sierpinski Gasket; a few years later, this   
definition was extended to include Laplacians on a large class of self--similar fractal sets \cite{Kigami93}, 
known as {\it post critically finite} sets (p.c.f. sets). 
The method of {\it spectral decimation} introduced by Fukushima and Shima in the 1990's provides a way to 
evaluate the eigenvalues of Kigami's Laplacian. In general terms, this method consists in finding the eigenvalues of the self--similar fractal set 
by taking limits of eigenvalues of discrete Laplacians that act on some graphs that approximate the fractal. The spectral decimation method was applied 
in \cite{Fukushima-Shima} to the Sierpinski Gasket, in order to give an explicit construction which allows one to obtain the set of eigenvalues.
 In \cite{Shima96} it was shown that it is possible to apply the spectral 
decimation method to a large collection of p.c.f. sets, including the family of fractals known as {\it nested fractals} that was introduced 
by T. Lindstr\o m in \cite{Lindstrem}. In addition to the Sierpinski Gasket, the spectral decimation method has been applied 
in several specific cases of p.c.f. fractals (e.g. \cite{CSW,Drenning-Strichartz,Ford-Steinhurst,Metz2004,Zhou1}); also, the method has been proved useful 
to study the spectrum of particular fractals that are not p.c.f. (e.g. \cite{Bajorin1,Bajorin2}), and of fractafolds modeled on the Sierpinski Gasket 
\cite{Strichartz2003}. The spectral decimation method has also shown to be a very useful tool for the analysis of the structure of the spectra of Laplacians
of some fractals (e.g. \cite{Hare-Zhou,Strichartz2005,Zhou2}).
\medskip 

\begin{figure}
\begin{center}
  \includegraphics[scale=0.2]{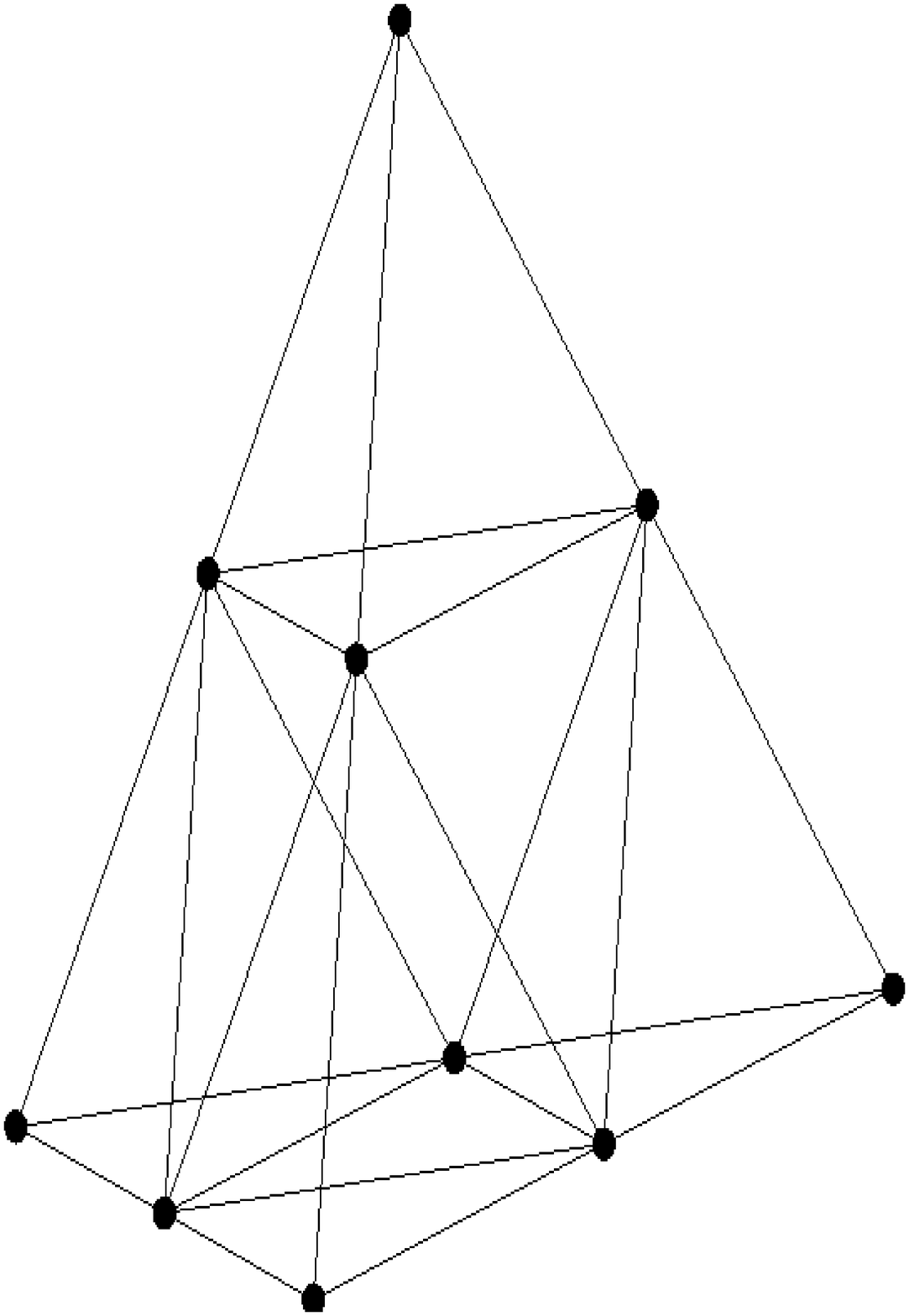}
\end{center}
  \caption{The first approximating graph $\Gamma_1$ for the fractal $\P 4$.}
 \label{ejemplo}
\end{figure}

\begin{figure}
\begin{center}
  \includegraphics[scale=0.2]{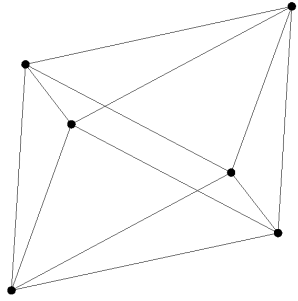}
\end{center}
  \caption{The graph $\Gamma_1$ for the fractal $\P 4$, minus the boundary points.}
 \label{sin_frontera}
\end{figure}

In the present work, we develope in an explicit way the spectral decimation method for an infinite collection of self--similar sets, that we will
denote by $\P n$ ($n\geq 2$ a positive integer). The definition of these sets is given in Definition \ref{conjuntos}. For the cases 
$n=2,3$, they correspond, respectively, to the unit interval and the Sierpinski Gasket. For larger values of $n$ they give a quite natural
extension of the Sierpinski Gasket to higher dimensions.  The spectral 
decimation method for the cases $n=2,3$ is presented with thorough detail in \cite{Strichartz}. 
Our presentation follows this reference to some extent. However, 
some technical difficulties arise for $n\geq 4$. This is mainly due to the fact that --even though the fractals considered are very symmetric --
the graphs approximating the fractal are not as homogeneous as the ones approximating the Sierpinski Gasket. For instance, if we consider the graph
obtained by by taking away the boundary points from $\Gamma_1$ (see Definition \ref{gammas} and Figures \ref{ejemplo} and \ref{sin_frontera}), then it will be a complete graph 
only for $n\leq 3$. A consequence of this, is the appereance of sets of two types of vertices that have to be dealt with separately, 
and which we denote by $F_{r,s}$ and $G_{r,s}$; for $n\leq 3$ the sets $G_{r,s}$ are empty. 
We also make the observation that the approximating graphs $\Gamma_k$ are non-planar when $n>3$.
\medskip 

In Section \ref{notacion}, we present general facts about self--similar sets, for the sake of completeness and in order to establish notation. 
In Section \ref{fractales} we introduce the sets $\P n$ that are the subject of study in this work; at the end of the section we find the Hausdorff dimension
of these fractals when embedded in Euclidean space. 
In Section \ref{graficas} we define the graphs that approximate
the self--similar sets $\P n$ and fix more notation. Our main result is presented in Section \ref{diseminar} (Theorem \ref{main}); it is shown that the 
eigenvalues and eigenfunctions of the discrete Laplacians of the approximating graphs can be obtained recursively. Finally, in Section \ref{laplacian}, 
it is shown that the eigenvalues of the Laplace operator in $\P n$ can be recovered by taking limits of the discrete Laplacians; 
in order to do this, we solve the so-called {\it renormalization problem} for this case (see Theorem \ref{renormalization}).

\section{Notation and Preliminaries}\label{notacion}

We denote by $\S n$ the shift--space with $n$ symbols. In this work we will always consider these $n$ symbols to be
the numbers $0,1,\dots,n-1$. $\S n$ is a compact space (see e.g. \cite{Kigami}) when equipped with the metric
$$
\delta(a_0 a_1 a_2 a_3\dots;b_0 b_1 b_2 b_3\dots)=r^k\qquad 0<r<1 
$$
where
$$
k=\min\{j\geq 0\ |\ a_j\neq b_j\}.
$$
We will use the dot notation $\dot a$, meaning that the symbol $a$ repeats to infinity.
\medskip 

Let ${\bf x}=x_0 x_1 x_2\dots$ an element of $\S n$, and $a\in\{0,\dots,n-1\}$.
We denote by $T_a$ the shift--operator given by 
$$
T_a({\bf x})=ax_0 x_1 x_2\dots
$$
It is easy to see that
$$
\delta\bigl(T({\bf x});T({\bf y})\bigr)=r\delta({\bf x};{\bf y}),
$$
so that $T_a$ is a contraction (by factor $r$). The space $\S n$ is a self--similar set, equal to $n$ smaller copies of itself, 
with $\{T_0,\dots,T_{n-1}\}$ the corresponding contractions. 
Even more so, it can be proved (Theorem 1.2.3. in \cite{Kigami}) that if $K$ is any self--similar set then it is homeomorphic 
to a quotient space of the form $\S n/\sim$ for a suitable equivalence relation.
\medskip 

For $K=\S n/\sim$, and ${\bf a}$ a word of length $m$ 
$$
{\bf a}=(a_0 a_1\dots a_{m-1})
$$ 
denote by $T_{\bf a}$ the shift--operator given by
$$
T_{\bf a}({\bf x})=a_0 \dots a_{m-1}x_0 x_1 x_2\dots
$$
The operator $T_{\bf a}$ is called an $m$-{\it contraction}, and the sets of the form $T_{\bf a}(K)$ are known as the {\it cells of level} $m$
of the self--similar set $K$. We note that, for each choice of $m$, $K$ is the union of the $n^m$ cells of level $m$.
\medskip 

\section{The Self--Similar Fractals $\P n$}\label{fractales}

Here we will introduce the self--similar fractals $\P n$ that are the subject of analysis in this work.
\medskip 

\begin{definicion}\label{conjuntos}
 For $n\in\N$ define $\P n$ as the quotient space $\S n/\sim$, with the equivalence relation given by
 $$
  a_0 a_1 a_2\dots a_k b\dot c\sim a_0 a_1 a_2\dots a_k c\dot b
 $$
for any choice of symbols $a_j$, $b$ and $c$.
 \end{definicion}

 $\P 1$ is a trivial space with only one element, $\P 2$ is homeomorphic to a compact interval in $\R$, and $\P 3$ is homeomorphic
 to the well known Sierpinski Gasket. For any value of $n$, $\P n$ can be embedded in Euclidean space; more precisely, there exists a (quite natural) 
 homeomorphism between $\P n$ and a compact self--similar set $K_n\subset R^{n-1}$. Below we define the sets $K_n$; 
 for these representations of $\P n$, we will be able to find their Hausdorff dimensions. 
\medskip

\begin{figure}
\begin{center}
\includegraphics[scale=0.25]{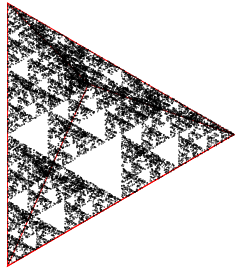}
\end{center}
 \caption{A representation of $K_4$, generated with MatLab.}
 \label{k4}
\end{figure}

Take $n$ points $x_0,\dots,x_{n-1}\in{\mathbb R}^{n-1}$ that do not lie in the same $(n-2)-$ dimensional hyperplane; for $n=3$
those points will be the ``vertices'' of the Sierpinski Gasket. For $n=4$ the fractal $K_4$ will be some sort of Sierpinski tetrahedron 
(see Figure \ref{k4}), while the four points $x_j$ will be the vertices of the tetrahedron.
\medskip 

Consider the contractions
$$
f_i(x)=\frac {x+x_i} 2\qquad i=1,\dots,{n-1}.
$$
We note that $f_i$ maps each $x\in\R^{n-1}$ to the midpoint of $x$ and $x_i$ (hence, leaving $x_i$ fixed). Define $K_n$ as the unique compact set
such that
$$
K_n=\bigcup_{i=0}^{n-1} f_i(K_n).
$$
We note that, for $i\neq j$, the sets $f_i(K_n)$ and $f_j(K_n)$ intersect at exactly one point: $f_i(x_j)=f_j(x_i)$. From this, it follows that
the map $\pi:\P n\rightarrow K_n$ given by
$$ 
\pi(\omega_0\omega_1\omega_2\dots)=\bigcap_{m\geq 0} f_{\omega_0}\circ f_{\omega_1}\circ\dots\circ f_{\omega_m} (K_n) 
$$
is a well defined homeomorphism; also, for every $k=0,\dots, n-1$, the following diagram commutes (cf. Theorem 1.2.3 in \cite{Kigami}): 
$$
\begin{CD}
\P n @>T_k>> \P n\\
@V\pi VV @VV\pi V\\
K_n @>>f_k> K_n
\end{CD}
$$
\medskip 

The sets $K_n$ satisfy the Moran--Hutchinson {\it open set condition}; namely, that there exists a bounded non--empty open set 
$O\subset\P n$ such that
$$
f_i(O)\subset O,\quad \forall i\in\{0,\dots n-1\}
$$
and
$$
f_i(O)\cap f_j(O),\quad \forall i\neq j.
$$
Just take $O=K_n\setminus\{x_0,\dots,x_{n-1}\}$. From this and the fact that $K_n$ is equal to $n$ contractions of itself (by factor $1/2$), 
it follows from Moran's theorem (Corollary 1.5.9 in \cite{Kigami}) that the Hausdorff dimension of $K_n$, respect to Euclidean metric, 
is equal to $\log n/\log 2$.
\medskip 

We end this section with two relevant notes:
\begin{itemize}
 \item 
For some values of $n$, it might be possible to embed $K_n$ isometrically into Euclidean space of a dimension $m$ smaller than $n-1$. Of course,
the dimension of the fractal gives a restriction to the minimal value of $m$.
\item 
The representations $K_n$ are somehow useful to visualize the self--similar fractals $\P n$. However, this representation and its metric do not
play any role in the analysis carried out in the next sections; we will therefore will focus in the more abstract definition of $\P n$ given 
at the beginning of this section.
\end{itemize}
\bigskip 

\section{Graph approximations of Self-Similar Sets}\label{graficas}

 In this and the next sections, we consider the self--similar set $\P n$ defined above, for an arbitrary but 
fixed value of $n\geq 2$. 
\medskip 

Let $V_0$ be the set of points in $\P n$ that have the form $\dot k$ with 
$k=1,\dots, n-1$. We call $V_0$ the {\it boundary} of $\P n$. Likewise, for $m\in\N$ let $V_m$ be the subset of $\P n$ of points of the 
form $a_0\dots a_{m-1}\dot k$. In other words, $x\in V_m$ if and only if it belongs to the image of $V_0$ under some $m$-contraction. 
\medskip 

Next, we define the graphs that will approximate $\P n$.

\begin{definicion}\label{gammas}
Denote by $\Gamma_0$ the complete graph of $n$ vertices, with $V_0$ its set of vertices.  
For $m\in\N$, let $\Gamma_m$ be the graph with set of vertices 
$V_m$ and edge relation established by requiring $x$ to be connected with $y$ if and only if there exists an $m$-contraction $T_{\bf a}$ 
such that both points $x$ and $y$ are in $T_{\bf a}(V_0)$.
\end{definicion}

We can see that an equivalent formulation is that two vertices $x$ and $y$ share an edge in $\Gamma_m$ only when their first $m$ symbols coincide.
It is worth noting that even though $V_0\subset V_1\subset V_2\cdots$, the edge relation is never preserved; this follows from the fact 
that if $x\neq y$ are connected in $\Gamma_m$, then their $(m+1)$-th symbols cannot be equal, so that they will not be connected in $\Gamma_{m+1}$. 
\medskip 

For each $m\in\N$ let $\Delta_m$ be the graph--Laplacian on $\Gamma_m$. We consider the Laplacian as acting on a space with boundary. 
More precisely, for a real--valued function $u$ defined on $V_m$ and $x$ in $V_m\setminus V_0$: 
$$  
\Delta_m u(x)=\sum_{y\sim x} (u(x)-u(y)),
$$
with the sum over all vertices $y$ that share an edge with $x$; the boundary values remain unchanged. Also $u$ is an eigenfunction of $\Delta_0$ 
with  eigenvalue $\lambda$, if
$$
\Delta_m u(x)=\lambda u(x),\qquad \forall x\in V_m\setminus V_0.
$$
We denote by $E_m(\cdot,\cdot)$ the associated quadratic form (known as the {\it energy product} of the graph):
$$
E_m(u,v)=(\Delta_m u,v)=\sum_{x\sim y} (u(x)-u(y))(v(x)-v(y))
$$
for $u$ and $v$ real--valued functions defined on $V_m$, and the sum being taken over the pairs of vertices $(x,y)$ that are connected to
each other. Also, we use the abbreviation $E(u)=E(u,u)$.

\section{Spectral Decimation}\label{diseminar}

Let $m>1$, and suppose $u$ is an eigenfunction of $\Delta_{m-1}$, with eigenvalue $\lambda_{m-1}$. We will show that it is always 
possible to extend this function to the domain $V_m$ so that it will be an eigenfunction of $\Delta_m$ (with not the same eigenvalue). 
In order to do this, we will derive necessary conditions for the extension to be an eigenfunction; in the process, 
it will become clear that those conditions are also sufficient.

Suppose that $u$ is an eigenfunction of $\Delta_m$ with eigenvalue $\lambda_m$; 
we aim to write the values of $u_m$ in $V_m\setminus V_{m-1}$ in terms of its values in $V_{m-1}$. 
Without loss of generality, we can restrict ourselves
to the set $V_m\cap T_{\bf a}(\P n)$ for a fixed $(m-1)$-contraction $T_{\bf a}$; this is because the vertices of $\Gamma_m$ that belong to the set
$(V_m\setminus V_{m-1})\cap T_{\bf a}(\P n)$ are not connected to any vertices outside the cell $T_{\bf a}(\P n)$.
Denote the elements of this set by
\begin{equation}
x_{b,c}={\bf a}b\dot c.\qquad b,c=0,\dots {n-1}.\label{def_puntos}
\end{equation}

\begin{figure}
\begin{center}
\includegraphics[scale=0.25]{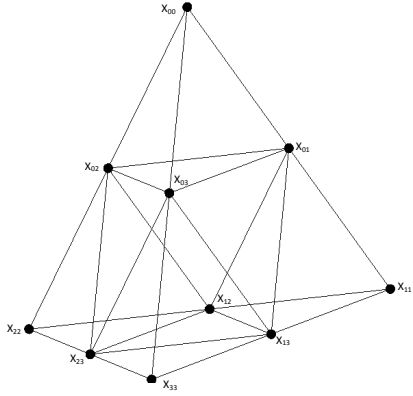}
\end{center}
 \caption{A cell of level $m-1$ of a graph $\Gamma_m$, approximating $\P 4$.}
 \label{notacion4}
\end{figure}

It is clear that $x_{b,c}=x_{c,b}$, and also that $x_{b,c}\in V_{m-1}$ if and only if $b=c$. This is shown in Figure \ref{notacion4}. 
\medskip

For each point $x_{r,s}\in V_m\cap T_{\bf a}(\P n)$ define the sets of vertices
\begin{align*}
 F_{r,s}&=\{x_{i,j}\ |\ i\neq j, \{r,s\}\cap\{i,j\}\neq\emptyset\}
 \\
 G_{r,s}&=\{x_{i,j}\ |\ i\neq j, \{r,s\}\cap\{i,j\}=\emptyset\}
\end{align*}
In other words, $F_{r,s}$ is the set of vertices (not in $\Gamma_{m-1}$) that are connected to the vertex $x_{r,s}$ in $\Gamma_m$, 
and $G_{r,s}$ is the set of vertices (not in $\Gamma_{m-1}$, either) that are not connected to it. 
\medskip 

\begin{figure}
\begin{center}
\includegraphics[scale=0.2]{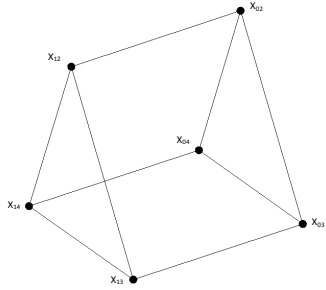}
\end{center}
 \caption{The graph determined by $F_{0,1}$ in a graph $\Gamma_m$, approximating $\P 5$.}
 \label{prismas}
\end{figure}

In the case $n=4$, the graph $\Gamma_m\setminus\Gamma_{m-1}$ is an octahedron; hence, for each pair $\{r,s\}$, the subgraph determined by the
vertices in $F_{r,s}$ is a $4$-cycle, while $G_{r,s}$ consists of a single vertex (the one opposite to $x_{r,s}$ in the octahedron). For general
$\P n$ we can see that:

\begin{itemize}
\item 
The graph $\Gamma_m\setminus\Gamma_{m-1}$ has $n(n-1)$ vertices, all of them with degree $2(n-2)$. Each of these vertices is connected to
another $2$ vertices in $\Gamma_{m-1}$.
 \item 
The subgraph determined by $F_{r,s}$ consists of two complete graphs, each one with $n-2$ vertices. The two complete graphs 
are joined to each other pairwise, thus forming a ``prism'', (a true prism only in the case $n=5$, where the base is a $3$-cycle, 
as shown in Figure \ref{prismas}).
 \item 
The subgraph determined by $G_{r,s}$ has $(n-2)(n-3)/2$ vertices, each one of them with degree $2(n-4)$. 
 \item 
In $\Gamma_m$, each vertex that belongs to $G_{r,s}$ is connected to exactly four vertices in $F_{r,s}$. On the other hand, each vertex
that belongs to $F_{r,s}$ is connected to $n-2$ vertices in $G_{r,s}$.
 \end{itemize}
\medskip

Now, having noted all that, we proceed with the calculations. For every $r\neq s$ we have 
\begin{equation}
 (2(n-1)-\lambda_m)u(x_{r,s})=u(x_{r,r})+u(x_{s,s})+\sum_{F_{r,s}} u(x_{i,j}).\label{Uno}
\end{equation}

Adding this up over all the possible values of $r$ and $s$, and rearranging terms yields

\begin{equation*}
 (2-\lambda_m)\sum_{r\neq s} u(x_{r,s}) =(n-1)\sum_{j=0}^{n-1} u(x_{j,j}),
\end{equation*}
which for any fixed $a\neq b$ can also be written in the form
\begin{equation*}
 (2-\lambda_m)\left(u(x_{a,b})+\sum_{F_{a,b}}u(x_{i,j})+\sum_{G_{a,b}}u(x_{i,j})\right)=(n-1)\sum_{j=0}^{n-1} u(x_{j,j}).
\end{equation*}
This, together with \eqref{Uno} allows us to express the sum of the values in $G_{a,b}$ in terms of $u(x_{a,b})$ and 
the values at points in $V_{m-1}$; namely, provided $\lambda_m\neq 2$, we have that
\begin{equation}  
\sum_{G_{a,b}}u(x_{i,j})=u(x_{a,a})+u(x_{b,b})
 -(2n-1-\lambda_m)u(x_{a,b})+\frac {(n-1)} {2-\lambda_m}\sum_{j=0}^{n-1} u(x_{j,j}).\label{Geab}
\end{equation}

Next, we will take the sum of the same terms, but only over the $x_{i,j}$ inside the set $F_{a,b}$ for fixed values $a\neq b$. 
Since $F_{a,b}$ contains two complete graphs with $n-1$ vertices and these complete graphs are pairwise connected to each other, 
it is clear that each $x_{i,j}\in F_{a,b}$ is connected to other $n-1$ vertices in $F_{a,b}$. Also, recall that each $x_{i,j}\in G_{a,b}$
is connected to exactly four vertices in $F_{a,b}$. For the vertices in $V_{m-1}$ we note that $x_{j,j}$ is connected to $n-2$ vertices 
in $x_{i,j}\in F_{a,b}$ if $j=a,b$ and to only two vertices otherwise. 
\medskip 

From the preceeding discussion it follows the equality 

\begin{align}
 (n-\lambda_m)\sum_{F_{a,b}}u(x_{i,j})&=4\sum_{G_{a,b}}u(x_{i,j}) + (n-2)(u(x_{a,a})+u(x_{b,b}))\notag \\
 &+2(n-2)u(x_{a,b})+2\sum_{j\neq a,b}u(x_{j,j})\label{Tres}
\end{align}

Consider the expresion given by \eqref{Uno} for $\{a,b\}=\{r,s\}$, multiply it by $n-\lambda_m$, and substitute equality \eqref{Tres} into it; 
this gives after arranging terms

\begin{align*}
 (\lambda_m^2-(3n-2)\lambda_m+2(n^2-2n+2))u(x_{a,b})&=4\sum_{G_{a,b}}u(x_{i,j})\\
 +2\sum_{j\neq a,b}u(x_{j,j})+(2(n-1)-&\lambda_m)(u(x_{a,a})+u(x_{b,b})).
\end{align*}
We want to get rid of the terms corresponding to $G_{a,b}$, so we replace it by \eqref{Geab}.
After straightforward computations, we can see that for $\lambda_m\neq 2$
\begin{align*}
(\lambda_m^2-(3n&+2)\lambda_m+2n(n+2))u(x_{a,b})=\\
&=\frac {\left(\lambda_m^2-2(n+2)+8n\right)(u(x_{a,a})+u(x_{b,b}))
+2(2n-\lambda_m)\sum_{j\neq a,b}u(x_{j,j})}{2-\lambda_m}.
\end{align*}
The quadratic equation for $\lambda_m$ in the left hand side has roots $n+2$ and $2n$. The one in the right hand side has roots $4$ and $2n$.
This gives us the following expression for $u(x_{a,b})$ in terms of the values of $u$ in $V_{m-1}$:
\begin{equation}
 u(x_{r,s})=\frac {(4-\lambda_m)(u(x_{r,r})+u(x_{s,s}))+2\sum_{j\neq r,s}u(x_{j,j})}
 {(2-\lambda_m)((n+2)-\lambda_m)}\label{finale}
\end{equation}
valid for any eigenvalue $\lambda_m\neq 2, n+2, 2n$. 
\medskip 

For $\lambda_m=0$ this reduces to
\begin{equation}
 u(x_{r,s})=\frac 2 {n+2}(u(x_{r,r})+u(x_{s,s}))+\frac 1 {n+2}\sum_{j\neq r,s}u(x_{j,j}).\label{caso_cero}
\end{equation}

It is clear from the construction that if $u(x_{a,b})$ is defined 
by \eqref{finale} and $\lambda_m$ is given by \eqref{recursiva}, then we have that
$$
\Delta_m u(x_{a,b})=\lambda_m u(x_{a,b}).\qquad (a\neq b).
$$

It remains to verify that this is valid as well in $V_{m-1}$. Of course, this cannot be true for arbitrary values of $\lambda_m$, but
only at most for specific values depending on $\lambda_{m-1}$; we will find those values in what follows. 
\medskip 

Take a point in $V_{m-1}$, say
$$
x_{p,p}={\bf a}\dot p,\qquad {\bf a}=a_0\dots a_{m-1}.
$$
Suppose that $a_k=q$ is the last symbol in ${\bf a}$ that is different from $p$;  we can assume that such symbol exists, since otherwise 
$x_{p,p}$ would be in the boundary $V_0$. With this, the point $x_{p,p}$ can also be written in the form
$$
x_{p,p}={\bf a'}\dot q,\qquad {\bf a'}=a_0\dots a_{k-1}pq\dots q
$$
with the necessary number of $q$'s to make ${\bf a'}$ a word of length $m-1$. Hence,  $x_{p,p}$ is in exactly two 
different $(m-1)$-cells: $T_{\bf a}(\P n)$ and $T_{\bf a'}(\P n)$, corresponding to each one of its two representations.   

\begin{figure}
\begin{center}
 \includegraphics[scale=0.3]{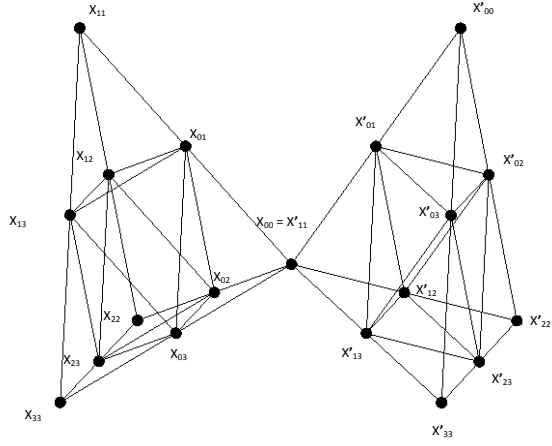}
\end{center}
 \caption{Two cells of level $m-1$ intersecting in a vertex of $\Gamma_m$.}
 \label{juntos}
\end{figure}

Denote by $x_{r,s}'$ the points in $T_{\bf a'}(\P n)\cap V_m$, defined as in \eqref{def_puntos} for the
points in $T_{\bf a}(\P n)\cap V_m$; in particular $x_{p,p}=x_{q,q}'$ (see Figure \ref{juntos}). The value of $u$ in the points  $x_{r,s}'$ is
given by the analogue of equation \eqref{finale}. The vertex $x_{p,p}$ is connected in $\Gamma_m$ to the $2(n-1)$ points of the form 
$x_{p,j}$ and $x_{q,j}'$, from which it follows that $u$ is an eigenfunction of $\Delta_m$ with eigenvalue $\lambda_m$, if and only if 
\eqref {finale} holds for all $x_{r,s}\in V_m\setminus V_{m-1}$ and the following equality holds for all $x_{p,p}\in V_{m-1}$:
\begin{equation}
 (2(n-1)-\lambda_m)u(x_{p,p})=\sum_{i\neq p, j\neq q} \left(u(x_{p,j})+u(x_{q,j}')\right).\label{Laplacian_m}
\end{equation}

On the other hand, since we know that $u$ is an eigenfunction of $\Delta_{m-1}$ with eigenvalue $\lambda_{m-1}$, 
we also have that:
\begin{equation}
 (2(n-1)-\lambda_{m-1}) u(x_{p,p})=\sum_{i\neq p,j\neq q} \left(u(x_{i,i})+u(x_{j,j}'\right).\label{Laplacian_m1}
\end{equation}

Replacing each term in the right hand side of \eqref{Laplacian_m} by its expression given by \eqref{finale} we can see that
\begin{equation*}
 (2(n-1)-\lambda_m)u(x_{p,p})=\frac {2(n-1)(4-\lambda_m)u(x_{p,p})+(2n-\lambda_m)}{(2-\lambda_m)((n+2)-\lambda_m)},
\end{equation*}
and using \eqref{Laplacian_m1} this gives
\begin{equation*}
 u(x_{p,p})=\frac{\left[2(n-1)(4-\lambda_m)+(2n-\lambda_m)(2(n-1)-\lambda_m)\right] u(x_{p,p})}
 {(2(n-1)-\lambda_m)(2-\lambda_m)((n+2)-\lambda_m)} 
\end{equation*}

Taking $u(x_{p,p})\neq 0$, and cancelling out, after computations the above equality reduces to the quadratic
\begin{equation}
 \lambda_m^2-(n+2)\lambda_m+\lambda_{m-1}=0,\label{cuadratic}
\end{equation}
which in turn gives the following recursive characterization of the eigenvalues:
\begin{equation}
 \lambda_m=\frac{(n+2)\pm\sqrt{(n+2)^2-4\lambda_{m-1}}} 2\label{recursiva}
\end{equation}

Since this procedure can be reversed, we have proved the following result.

\begin{teorema}\label{main}
 Let $\lambda_m\neq 2, n+2, 2n$, and let $\lambda_{m-1}$ be given by \eqref{cuadratic}. 
 Suppose $u$ is an eigenfunction of $\Delta_{m-1}$ with eigenvalue $\lambda_{m-1}$. Extend $u$ to $V_m$ by \eqref{finale}.
 Then $u$ is an eigenfunction of $\Delta_m$ with eigenvalue $\lambda_m$. Conversely, if $u$ is an eigenfunction of 
 $\Delta_m$ with eigenvalue $\lambda_m\neq 2, n+2, 2n$, then its restriction to $V_{m-1}$ is
 an eigenfunction of $\Delta_{m-1}$ with eigenvalue $\lambda_{m-1}$.
 \end{teorema}

\section{The Laplacian on the Self--Similar Fractals}\label{laplacian}

In order to define the Laplace operator of a p.c.f. fractal by means of graph approximations, it is required
to solve the  so called  {\it renormalization problem} for the fractal (e.g. \cite{Strichartz}, Chapter 4); 
roughly, this consists in normalizing the graph energies in $\Gamma_m$ in order to obtain a self--similar energy in the fractal by taking the limit. 
This can be achieved if the energies are such that they remain constant for 
each harmonic extension from $\Gamma_m$ to $\Gamma_{m+1}$. Below, we do this for the $\P n$ sets.

\begin{definicion}
 For a given function $u$ with domain $V_{m-1}$, we call the extension of $u$ to $V_m$ given by \eqref{caso_cero} its {\it harmonic extension}.
\end{definicion} 

The next result, gives the explicit solution of the renormalization problem the for $\P n$. 

\begin{teorema}\label{renormalization}
 Let $u:V_{m-1}\rightarrow\R$ arbitrary, and let $u':V_m\rightarrow\R$ be its harmonic extension. Then
 $$ 
 E_m(u')=\frac n {n+2} E_{m-1}(u)
 $$
\end{teorema}
{\it Proof.} 
Note that the energy at level $k$ of a given function equals the sum of the energies at all the $k'$--cells for any $k'\leq k$,
since different cells share no edges. This allows to restrict ourselves to one fixed $m-1$--cell both while considering $E_m(u')$ 
and $E_{m-1}(u)$. We use the notation of the previous section for the vertices of $\Gamma_m$ in this cell, and write $\tilde E$
for the energy restricted to this cell. We can readily see that 
\begin{align}
 \tilde E_{m-1}(u)&=\sum_{i\neq j}(u(x_{i,i})-u(x_{j,j}))^2.\notag\\ 
&=(n-1)\sum_{i=0}^{n-1}u^2(x_{i,i})-2\sum_{i\neq j}u(x_{i,i}u(x_i,j))\label{anterior}.
 \end{align} 
In order to evaluate the energy $\tilde E_m$ we consider first the edges joining vertices in $V_{m-1}$ with vertices in $V_m\setminus V_{m-1}$:
The edge joining the vertex $x_{a,a}$ with the vertex $x_{a,k}$ contributes to the energy by
\begin{equation*}
(u(x_{a,a})-u(x_{a,k}))^2=\frac 1 {(n+2)^2}\left(n u(x_{a,a})-2 u(x_{k,k})-\sum_{j\neq a,k} u(x_{j,j})\right)^2. 
\end{equation*}
When adding up this over all possible pairs $a\neq k$, each $x_{r,r}$ will appear $n-1$ times as the $x_{a,a}$, another $n-1$ times as the $x_{k,k}$ and
$(n-1)(n-2)$ as one of the $x_{j,j}$'s. Each pair double product $2x_{r,r}x_{s,s}$ will appear twice for $\{r,s\}=\{a,k\}$, $2(n-2)$ times for
$\{r,s\}=\{a,j\}$ for some $j$, also $2(n-2)$ times for $\{r,s\}=\{k,j\}$ for some $j$, and finally $(n-2)(n-3)$ times 
when both $r$ and $s$ are one of the $j$'s. All this implies that the contribution to the energy from these edges is, after simplification:
\begin{align}
\sum_{a\neq k}(u(x_{a,a})&-u(x_{a,k}))^2\notag\\
&=\frac {(n-1)(n^2+n+2)} {(n+2)^2}\sum_{i=0}^{n-1} u^2(x_{i,i})-\frac{2(n^2+n+2)}{(n+2)^2}\sum_{i\neq j}x_{i,i}x_{j,j}.\label{energia1}
\end{align}
On the other hand, the contribution from the edge that joins the vertices $x_{a,b}$ and $x_{a,c}$ (in $V_m\setminus V_{m-1}$) equals
\begin{equation*}
(u(x_{a,b})-u(x_{a,c}))^2=\frac 1 {(n+2)^2}\left(u(x_{b,b})-u(x_{c,c})\right)^2. 
\end{equation*}
Taking the sum over all of the vertices in $V_m\setminus V_{m-1}$ yields
\begin{align}
\sum_{a\neq b\neq c} &(u(x_{a,b})-u(x_{a,c}))^2\notag\\
&=\frac {(n-1)(n-2)} {(n+2)^2}\sum_{i=0}^{n-1} u^2(x_{i,i})-\frac {2(n-2)} {(n+2)^2}\sum_{i\neq j}u(x_{i,i})u(x_{j,j}).\label{energia2}
\end{align}
From \eqref{energia1} and \eqref{energia2} it follows that the total energy of the cell is
\begin{equation}
\tilde E_m(u')=\frac {n(n-1)}{n+2}\sum_{i=0}^{n-1} u^2(x_{i,i})-2n(n+2).
\end{equation}
This, together with \eqref{anterior} gives 
$$ 
\tilde E_m(u')=\frac n {n+2} \tilde E_{m-1}(u).
$$
Taking this result for all the $m-1$-cells concludes the proof.

\rightline{$\square$}

\begin{definicion}
 The energy in $\P n$ is given by
 $$
 E(u)=\lim_{m\to\infty}\left(\frac {n+2} n\right)^m E_m(u).  
 $$
 The domain of $E(\cdot)$ being the space $D^{(n)}$ of functions such that the energy is finite. Write $D^{(n)}_0$ for the subspace of 
 $D^{(n)}$ of functions that vanish on the boundary. The energy product $E(u,v)$ can be recovered by the polarization identity. 
\end{definicion}

Let $\mu$ be a self--similar measure in $\P n$, the Laplacian $\Delta_\mu$ is given by: 

\begin{definicion}{\bf (Kigami's Laplacian)}
With $\mu$ and $\Delta_\mu$ as above, we say that $u$ is in the domain of $\Delta_\mu$ if there exsists a continuous function $f$ such that
$$ 
E(u,v)=-\int_{P_n} fv\ d\mu\qquad\forall v\in D^{(n)}_0.
$$
In such case, we define $\Delta_\mu u=f$. 
\end{definicion}

Aside from the above weak representation, a pointwise formula can be obtained for $\Delta_mu$, proceeding in exactly in the same
way as in \cite{Strichartz} (Theorem 2.2.1). In the case where $\mu$ is the standard measure in $P_n$ (i.e. the only Borel regular measure such that
the measure of every $m$-cell is equal to $n^{-m}$), the pointwise formula is
$$
\Delta_\mu u(x)=\frac n 2 \lim_{m\to\infty}(n+2)^m \Delta_m u(x).
$$

This leads to the following: If a sequence $\{\lambda_m\}$ is defined recursively by \eqref{recursiva}
(assuming that $\lambda_m$ is never equal to $n, n+2$ or $2n$), and $u_m$ is given by relation \eqref{finale} then 
\begin{equation*}
\lambda=\frac n 2 \lim_{m\to\infty}(n+2)^m \lambda_m 
\end{equation*}
is an eigenvalue of $\Delta_\mu$ with eigenfunction $u$ given by the limit $u_m\to u$. The limit above exists provided that the sign 
in relation \eqref{recursiva} is chosen to be ``+'' for at most a finite number of times.
\bigskip 

{\bf Acknowledgements.} We are very grateful to Alejandro Butanda and Yolanda Ortega for their valuable assistance with the figures, and to
the referee for many useful comments.

\end{document}